\documentstyle[12pt,fullpage]{amsproc}  

\renewcommand{\mod}[1]{{\ifmmode\text{\rm\ (mod~$#1$)}\else\discretionary{}{}{\hbox{ }}\rm(mod~$#1$)\fi}}

\newsymbol\dnd 232D

\newcommand{\V}{\mathop{\rm Var}}
\renewcommand{\Re}{\mathop{\rm Re}}

\begin{document}

\title{Asymmetries in the Shanks--R\'enyi Prime~Number~Race}
\author{Greg Martin}
\address{Department of Mathematics\\University of Toronto\\Canada M5S 3G3}
\email{gerg@@math.toronto.edu}
\begin{abstract}
It has been well-observed that an inequality of the type $\pi(x;q,a) >
\pi(x;q,b)$ is more likely to hold if $a$ is a non-square modulo $q$ and $b$
is a square modulo $q$ (the so-called ``Chebyshev Bias''). For instance, each
of $\pi(x;8,3)$, $\pi(x;8,5)$, and $\pi(x;8,7)$ tends to be somewhat larger
than $\pi(x;8,1)$. However, it has come to light that the tendencies of these
three $\pi(x;8,a)$ to dominate $\pi(x;8,1)$ have different strengths. A
related phenomenon is that the six possible inequalities of the form
$\pi(x;8,a_1) > \pi(x;8,a_2) > \pi(x;8,a_3)$ with $\{a_1,a_2,a_3\}=\{3,5,7\}$
are not all equally likely---some orderings are preferred over others. In
this paper we discuss these phenomena, focusing on the moduli $q=8$ and
$q=12$, and we explain why the observed asymmetries (as opposed to other
possible asymmetries) occur.
\end{abstract}
\maketitle

\section{Background}

Let $\pi(x;q,a)$ denote the number of primes not exceeding $x$ that
are congruent to $a$ modulo $q$. We know from the prime number theorem
in arithmetic progressions that the two counting functions
$\pi(x;q,a)$ and $\pi(x;q,b)$ are asymptotically equal as $x$ tends to
infinity (as long as $a$ and $b$ are both coprime to $q$).  However,
more complicated behavior emerges when we compare these counting
functions for finite values of $x$. Imagine $\pi(x;q,a)$ and
$\pi(x;q,b)$ as representing the two contestants in a race; as the
primes are listed in order, the contestant $\pi(x;q,a)$ takes a step
each time a prime is congruent to $a$ mod $q$, and similarly for
$\pi(x;q,b)$. How often is the first contestant ahead of the second?
This ``race game'' is easily extended to include several contestants.

As a prime example (!), consider the two contestants $\pi(x;4,1)$ and
$\pi(x;4,3)$. (We won't pay any attention to the other contestant $\pi(x;4,2)$,
who, while quick out of the starting blocks, is rather lacking in endurance.)
Chebyshev was the first to note that there are ``many more'' primes congruent
to 3\mod 4 than to 1\mod 4. Indeed, the first value of $x$ for which
$\pi(x;4,1)>\pi(x;4,3)$ is $x=2$6,861 (Leech~\cite{leech}). Even this victory
for $\pi(x;4,1)$ is short-lived, as 26,861 is the first of a pair of twin
primes, and so $\pi(x;4,3)$ catches right back up and does not relinquish the
lead again until
$x=6$16,481.

Similar biases are observed in race games to other moduli, especially to
small moduli. For example, $\pi(x;3,1)$ does not exceed $\pi(x;3,2)$ for the
first time until $x=6$08,981,813,029 (Bays and Hudson~\cite{BH3}). We can
also compare the four counting functions $\pi(x;8,a)$ for $a\in\{1,3,5,7\}$.
By the time $x$ equals 271, each of $\pi(x;8,3)$, $\pi(x;8,5)$, and
$\pi(x;8,7)$ has been in first place; but $\pi(x;8,1)$ does not even obtain
undisputed possession of third place in this four-way race until
$x=5$88,067,889 (Bays and Hudson~\cite{BH8}).

All of these biases just mentioned are instances of a universal
tendency for contestants $\pi(x;q,a)$ where $a$ is not a square modulo
$q$ to run ahead of contestants $\pi(x;q,b)$ where $b$ is a square
modulo $q$. We briefly indicate why this is the case through an
analytic argument (though see Hudson \cite{common} for a different
approach), at the same time establishing some of the notation to be
used throughout this paper. We will always assume that the modulus $q$
is fixed, and therefore we will not care about the dependence of
implicit $O$-constants on $q$.

Let $\psi(x;q,a)$ have its usual meaning,
\begin{equation*}
\psi(x;q,a) = \sum\begin{Sb}n\le x \\ n\equiv a\mod q\end{Sb} \Lambda(n) = 
\sum\begin{Sb}p^r\le x \\ p^r\equiv a\mod q\end{Sb} \log p,
\end{equation*}
and set $\psi(x)=\psi(x;1,1)$. Under the Generalized Riemann
Hypothesis (GRH), the explicit formula from the proof of the prime
number theorem for arithmetic progressions (see~\cite{davenport})
gives
\begin{equation}
\psi(x;q,a) = {\psi(x)\over\phi(q)} - {1\over\phi(q)} \sum
\begin{Sb}\chi\mod q \\ \chi\ne\chi_0\end{Sb} \bar\chi(a)
\sum_\gamma {x^{1/2+i\gamma}\over1/2+i\gamma} + O(\log^2 x)
\label{psiformula}
\end{equation}
as long as $a$ is coprime to $q$. Here the inner sum is indexed by the
imaginary parts $\gamma$ of the nontrivial zeros of the Dirichlet
$L$-function corresponding to the character $\chi$, and should be
interpreted as
\begin{equation}
\lim_{T\to\infty} \sum_{|\gamma|<T} {x^{1/2+i\gamma}\over1/2+i\gamma}
\label{liminterp}
\end{equation}
so that it will converge.

We isolate the contribution to $\psi(x;q,a)$ from the primes themselves,
defining
\begin{equation*}
\theta(x;q,a) = \sum\begin{Sb}p\le x \\ p\equiv a\mod q\end{Sb} \log p,
\end{equation*}
so that
\begin{multline}
\theta(x;q,a) = \psi(x;q,a) - \sum_{b^2\equiv a\mod q} \theta(x^{1/2};q,b) \\
{}- \sum_{c^3\equiv a\mod q} \theta(x^{1/3};q,c) - \dots.
\label{stork}
\end{multline}
Notice that the number of terms in the first sum on the right-hand side of
this equation equals the number of square roots modulo $q$ possessed by $a$
(in particular, the sum is empty if $a$ is a non-square modulo $q$). Let us
define
\begin{equation}
c(q,a) = {-1} + \# \{ 1 \leq b \leq q\colon b^2 \equiv a \bmod q \},
\label{cqadef}
\end{equation}
so that $c(q,a)$ is an extension of the Legendre symbol $(\frac aq)$ to all
moduli $q$ (though not an extension with nice multiplicativity properties);
then the number of terms in the first sum on the right-hand side of
equation~(\ref{stork}) is exactly $c(q,a)+1$.

Invoking the prime number theorem for arithmetic progressions, we have
$\theta(y;q,a) = y/\phi(q) + O(\sqrt y\log^2y)$ for a fixed modulus
$q$ (assuming GRH). Using this fact together with the explicit
formula~(\ref{psiformula}), equation~(\ref{stork}) becomes
\begin{multline}
\theta(x;q,a) = {\psi(x)\over\phi(q)} - {1\over\phi(q)} \sum
\begin{Sb}\chi\mod q \\ \chi\ne\chi_0\end{Sb} \bar\chi(a) \sum_\gamma
{x^{1/2+i\gamma}\over1/2+i\gamma} \\
{}- (c(q,a)+1){\sqrt x\over\phi(q)} + O(x^{1/3}).
\label{thetaformula}
\end{multline}
In particular, we have $\theta(x) = \theta(x;1,1) = \psi(x)-\sqrt x +
O(x^{1/3})$.

Converting equation~(\ref{thetaformula}) to a formula for $\pi(x;q,a)$
involves only a straightforward partial summation argument. We phrase the
final result in terms of a normalized error term for $\pi(x;q,a)$, namely
\begin{equation}
E(x;q,a) = {\log x\over\sqrt x} \big( \phi(q)\pi(x;q,a)-\pi(x) \big).
\label{Exqadef}
\end{equation}
From equation~(\ref{thetaformula}) applied to both $\theta(x;q,a)$ and
$\theta(x)$, one can derive \cite[Lemma 2.1]{RS}
\begin{equation}
E(x;q,a) = -c(q,a) - \sum\begin{Sb}\chi\mod q \\ \chi\ne\chi_0\end{Sb}
\bar\chi(a)E(x,\chi) + O\big( {1\over\log x} \big),
\label{Exqaformula}
\end{equation}
where we have defined
\begin{equation}
E(x,\chi) = \sum_\gamma {x^{i\gamma}\over1/2+i\gamma}  \label{Exchidef}
\end{equation}
(interpreted similarly to~(\ref{liminterp}) to ensure its conditional
convergence). The behavior of $E(x;q,a)$ is therefore that of a
function oscillating in a roughly bounded fashion about the mean value
$-c(q,a)$, which is positive if $a$ is a non-square\mod q and negative
if $a$ is a square\mod q. These two different possible mean values are
the source of the bias towards nonsquare contestants.

Rubinstein and Sarnak \cite{RS} have quantified these biases. Define
$\delta_{q;a_1,\dots,a_r}$ to be the logarithmic density of the set of
real numbers $x$ such that the inequalities
\begin{equation*}
\pi(x;q,a_1) > \pi(x;q,a_2) > \dots > \pi(x;q,a_r)
\end{equation*}
hold, where the logarithmic density of a set $S$ is
\begin{equation}
\lim_{x\to\infty} \, {1\over\log x} \int_{[2,x]\cap S} {dt\over t}
\label{logdensitydef}
\end{equation}
assuming the limit exists. Let us assume not only GRH, but also that the
nonnegative imaginary parts of the nontrivial zeros of all Dirichlet
$L$-functions are linearly independent over the rationals, a hypothesis we
shall abbreviate LI. Under these assumptions, Rubinstein and Sarnak proved
that $\delta_{q;a_1,\dots,a_r}$ always exists and is strictly positive. They
also proved that $\delta_{q;a,b}>\frac12$ if and only if $a$ is a nonsquare\mod
q and $b$ is a square\mod q, and calculated several of these densities; for
example, $\delta_{4;3,1}=0.9959\dots$ and $\delta_{3;2,1}=0.9990\dots$.

In joint work with Feuerverger \cite{FM}, we calculated many other densities
(under the hypotheses GRH and LI). One of our discoveries was that the
numerical values of the densities can vary even when the modulus $q$ is
fixed. For example, modulo 8 the only square is 1 while the three
nonsquares are 3, 5, and 7, and modulo 12 the only square is 1 while the
three nonsquares are 5, 7, and 11; we calculated that
\begin{equation}
\text{
$
\begin{aligned}
\delta_{8;3,1} &= 0.999569 \\
\delta_{8;7,1} &= 0.998938 \\
\delta_{8;5,1} &= 0.997395
\end{aligned}
$
\qquad and\qquad
$
\begin{aligned}
\delta_{12;11,1} &= 0.999977 \\
\delta_{12;5,1} &= 0.999206 \\
\delta_{12;7,1} &= 0.998606.
\end{aligned}
$}
\label{2waydeltas}
\end{equation}
We also found that for race games involving more than two contestants,
certain orderings of the contestants are more likely than others even
if the residue classes involved are all squares or all nonsquares\mod
q, a situation that was foreshadowed in~\cite{RS}. For example, we
calculated that
\begin{equation}
\text{
$
\begin{aligned}
\delta_{8;3,5,7} &= \delta_{8;7,5,3} =  0.192801 \\
\delta_{8;3,7,5} &= \delta_{8;5,7,3} =  0.166426 \\
\delta_{8;5,3,7} &= \delta_{8;7,3,5} =  0.140772
\end{aligned}
$
\quad and\quad
$
\begin{aligned}
\delta_{12;5,7,11} &= \delta_{12;11,7,5} =  0.198452 \\
\delta_{12;7,5,11} &= \delta_{12;11,5,7} =  0.179985 \\
\delta_{12;5,11,7} &= \delta_{12;7,11,5} =  0.121563.
\end{aligned}
\qquad
$}
\label{3waydeltas}
\end{equation}

\noindent The goal of this paper is to begin to understand these
recently discovered types of asymmetries. We will focus on the densities
listed in equations~(\ref{2waydeltas}) and~(\ref{3waydeltas}),
explaining how we could have predicted that $\delta_{8;5,1}$ would be
smaller than both $\delta_{8;3,1}$ and $\delta_{8;7,1}$, for
example. The hope is that a very concrete inspection of these special
cases will function as a starting point for a future analysis of the
general case.

\section{Error terms and random variables}

The great utility of the hypothesis LI, concerning the linear independence of
the nonnegative imaginary parts of the nontrivial zeros of Dirichlet
$L$-functions, is in facilitating calculations involving those zeros that
are based on harmonic analysis. In this paper we will phrase these
calculations in terms of random variables, focusing on underscoring the ideas
involved rather than belaboring the analytic details.

Notice that we can write
\begin{equation*}
E(x,\chi) = \sum_{\gamma>0} \big( {x^{i\gamma}\over1/2+i\gamma} +
{x^{-i\gamma}\over1/2-i\gamma} \big) = 2\sum_{\gamma>0} {\sin(\gamma\log x +
\alpha_\gamma)\over\sqrt{1/4+\gamma^2}}
\end{equation*}
for certain real numbers $\alpha_\gamma$ independent of $x$. The hypothesis
LI implies that $1/2$ can never be a zero of $L(s,\chi)$, which is why we do
not have to consider $\gamma=0$. Also, under LI, any vector of the form
\begin{equation*}
\{\sin(\gamma_1\log x+\alpha_{\gamma_1}), \dots, \sin(\gamma_k\log
x+\alpha_{\gamma_k})\}
\end{equation*}
becomes uniformly distributed over the
$k$-dimensional torus as $x$ tends to infinity. (It is the presence of $\log
x$ in this statement that requires us to define the
$\delta_{q;a_1,\dots,a_r}$ as logarithmic densities rather than natural
densities.) It can be shown that $E(x,\chi)$ has a limiting (logarithmic)
distribution as $x$ tends to infinity; moreover, this distribution can also
be described in terms of random variables. We now give this description.

For any positive number $\beta$, let $Z_\beta$ be a random variable that is
uniformly distributed on the unit circle in the complex plane; we make the
convention that $Z_{-\beta} = \overline{Z_\beta}$. We stipulate that
any collection $\{Z_{\beta_i}\}$ with all of the $\beta_i$ distinct and positive
is an independent collection of random variables. For any Dirichlet character
$\chi\mod q$, define the random variable
\begin{equation*}
X(\chi) = \sum_\gamma {Z_\gamma\over\sqrt{1/4+\gamma^2}},
\end{equation*}
where again the sum is indexed by the imaginary parts of the nontrivial zeros
of $L(s,\chi)$. We can also write
\begin{equation}
X(\chi) = 2 \sum_{\gamma>0} {X_\gamma\over\sqrt{1/4+\gamma^2}} 
\label{Xchidef}
\end{equation}
(since $L(\frac12,\chi)\ne0$), where the $X_\gamma=\Re Z_\gamma$ are
independent random variables each distributed on $[-1,1]$ with the
sine distribution. One can then show (see \cite{RS}), assuming GRH and LI,
that {\it the limiting distribution of $E(x,\chi)$ is identical to the
distribution of the random variable~$X(\chi)$.} Similarly, it follows
from equation~(\ref{Exqaformula}) that the limiting distribution of
$E(x;q,a)$ is the same as the distribution of the random variable
\begin{equation*}
-c(q,a) + \sum\begin{Sb}\chi\mod q \\ \chi\ne\chi_0\end{Sb} X(\chi).
\end{equation*}
(One might expect the summand to be something like $\bar\chi(a) X(\chi)$
rather than simply $X(\chi)$, but the coefficient $\bar\chi(a)$ disappears
early in the argument because $\bar\chi(a)Z_\gamma$ is the same random
variable as $Z_\gamma$ itself.) We remark that the hypothesis LI
implies that the various $X(\chi)$ are mutually independent random variables.

Let us examine these normalized error terms and random variables more
concretely for the moduli $q=8$ and $q=12$. For a fundamental discriminant
$D$, let $\chi_D(n) = ({D\over n})$ using Kronecker's extension of the
Legendre symbol. Then the three nonprincipal characters\mod8 are
$\chi_{-8}$, $\chi_{-4}$, and $\chi_8$, while the three nonprincipal
characters\mod{12} are $\chi_{-4}$, $\chi_{-3}$, and $\chi_{12}$.
(Table~\ref{chartable} explicitly lists the values taken by these
characters. We shall abuse notation a bit and also denote by $\chi_D$ a
character modulo 8 or 12 that is induced by the primitive character $\chi_D$,
whose conductor is $|D|$.)

\begin{table}[b]
\begin{tabular}{|c|ccc|}
\hline
$\chi$ & $\chi(3)$ & $\chi(5)$ & $\chi(7)$ \\
\hline
$\chi_{-8}$ & 1 & $-1$ & $-1$ \\
$\chi_{-4}$ & $-1$ & 1 & $-1$ \\
$\chi_8$ & $-1$ & $-1$ & 1 \\
\hline
\end{tabular}
\hfil
\begin{tabular}{|c|ccc|}
\hline
$\chi$ & $\chi(5)$ & $\chi(7)$ & $\chi(11)$ \\
\hline
$\chi_{-4}$ & 1 & $-1$ & $-1$ \\
$\chi_{-3}$ & $-1$ & 1 & $-1$ \\
$\chi_{12}$ & $-1$ & $-1$ & 1 \\
\hline
\end{tabular}
\medskip
\caption{Values of the nonprincipal characters\mod8 and\mod{12}}
\label{chartable}
\end{table}

Now, if we want to consider how often $\pi(x;8,3)$ exceeds $\pi(x;8,1)$, for
example, we can look at the limiting distribution of the normalized
difference $E(x;8,3) - E(x;8,1)$ and ask what proportion of that distribution
lies above 0. From the explicit formula~(\ref{Exqaformula}), we have
\begin{equation*}
\begin{split}
\begin{split}
E(x;8,3) - E(x;8,1) &= 4 + \sum\begin{Sb}\chi\mod q \\ \chi\ne\chi_0\end{Sb}
(1-\bar\chi(3))E(x,\chi) + O\big( {1\over\log x} \big) \\
&= 4 + 2E(\chi_{-4}) + 2E(\chi_8) + O\big( {1\over\log x} \big).
\end{split}
\end{split}
\end{equation*}
Thus $E(x;8,3) - E(x;8,1)$ has the same limiting distribution as the
random variable $4+2X(\chi_{-4}) + 2X(\chi_8)$, where $X(\chi)$ is as
in equation~(\ref{Xchidef}). In particular, the density
$\delta_{8;3,1}$ equals the mass given to the interval $(0,\infty)$ by
this limiting distribution, or in other words simply
$\Pr(4+2X(\chi_{-4}) + 2X(\chi_8)>0)$. In fact, if we define
\begin{equation}
\begin{split}
X_{8;3,1} &= 4 + 2X(\chi_{-4}) + 2X(\chi_8) \\
X_{8;5,1} &= 4 + 2X(\chi_{-8}) + 2X(\chi_8) \\
X_{8;7,1} &= 4 + 2X(\chi_{-8}) + 2X(\chi_{-4})
\end{split}
\label{X8a1s}
\end{equation}
and
\begin{equation}
\begin{split}
X_{12;5,1} &= 4 + 2X(\chi_{-3}) + 2X(\chi_{12}) \\
X_{12;7,1} &= 4 + 2X(\chi_{-4}) + 2X(\chi_{12}) \\
X_{12;11,1} &= 4 + 2X(\chi_{-4}) + 2X(\chi_{-3}),
\end{split}
\label{X12a1s}
\end{equation}
then in each case, the distribution of the random variable
$X_{q;a,1}$ is the same as the limiting distribution of the difference
$E(x;q,a)-E(x;q,1)$, and $\delta_{q;a,1}=\Pr(X_{q;a,1}>0)$.

If we have several random variables, each with mean 4 and symmetric about
that mean, which ones will take positive values most often? If the random
variables have roughly the same shape, then we expect the ones with the
smallest variance to stay above 0 the most. So let's compute the variances of
the random variables $X_{q;a,1}$.

Any variance of the $\V(cX_\gamma)$ with $c>0$ is simply $\frac12c^2$,
and the various $X_\gamma$ are independent; so if we define
$V(\chi)=\V(X(\chi))$, we see from the definition~(\ref{Xchidef}) of
$X(\chi)$ that
\begin{equation}
V(\chi)=\sum_{\gamma>0} {2\over1/4+\gamma^2}.  \label{Vchidef}
\end{equation}
We know that the larger the conductor of a character is, the more
numerous and low-lying (close to the real axis) the zeros of the
corresponding $L(s,\chi)$ will be. In fact, the order of magnitude
of the sum in equation~(\ref{Vchidef}) is known to be the logarithm of
the conductor of $\chi$, at least on GRH; one can see this from the
formula (see Davenport \cite[p.~83]{davenport})
\begin{equation}
V(\chi) = \log{q\over\pi} - \gamma_0 - (1+\chi(-1))\log2 +
2\mathop{\rm Re}{L'(1,\chi)\over L(1,\chi)}
\label{davformula}
\end{equation}
when $\chi$ is a primitive character\mod q, where $\gamma_0$ is
Euler's constant. Therefore, $V(\chi)$ will be larger when the
conductor of $\chi$ is large. In particular, we should expect
\begin{equation}
V(\chi_{12}) > V(\chi_{-8}) > V(\chi_8) > V(\chi_{-4}) > V(\chi_{-3}),
\label{Vrelative}
\end{equation}
and the numerical computation of the variances verifies these
expectations (see Table~\ref{vartable}).

\begin{table}[t]
\begin{tabular}{|c|c|}
\hline
$\chi$ & $V(\chi)$ \\
\hline
$\chi_{-3}$ & $0.11323$ \\ 
$\chi_{-4}$ & $0.15557$ \\ 
$\chi_8$    & $0.23543$ \\ 
$\chi_{-8}$ & $0.31607$ \\ 
$\chi_{12}$ & $0.33017$ \\ 
\hline
\end{tabular}
\medskip
\caption{Values of $V(\chi)=\sum_{\gamma>0} {2\over1/4+\gamma^2}$}
\label{vartable}
\end{table}

Why do we say that we expect $V(\chi_{-8}) > V(\chi_8)$, when the two
characters have the same conductor? There is a secondary phenomenon,
namely that the zeros of $L$-functions corresponding to even
characters tend to be not as low-lying as those of $L$-functions
corresponding to odd characters (the trivial zero at $s=0$ of an
$L$-functions associated to an even character seems to have a
repelling effect on the nontrivial zeros). Indeed, the term $-
(1+\chi(-1))\log2$ in the formula~(\ref{davformula}) for $V(\chi)$,
which vanishes for odd characters $\chi$, slightly lowers the value of
$V(\chi)$ for even characters~$\chi$.

Of course this observation would be spurious if the behavior of the real part
of ${L'(1,\chi)/L(1,\chi)}$ were much different for odd and even
characters. While there is no reason to suspect that this should be the case,
it seems hard to say anything substantial about the distribution of these
values (this is a subject that warrants further investigation). Nevertheless,
a look at lists of the first several zeros of Dirichlet $L$-functions with
small conductor does confirm that the zeros of $L(s,\chi)$ are lower-lying
when $\chi$ is odd than when $\chi$ is even.

Returning to equations~(\ref{X8a1s}) and~(\ref{X12a1s}), we can easily compute
the variance of the random variables $X_{q;a,1}$ (again since the various
$X(\chi)$ are independent by LI). For example.
\begin{equation*}
\V(X_{12;5,1}) = 4V(\chi_{-3}) + 4\V(\chi_{12}) = 4W_{12} - 4V(\chi_{-4}),
\end{equation*}
where we define
\begin{equation*}
W_q = \sum\begin{Sb}\chi\mod q \\ \chi\ne\chi_0\end{Sb} V(\chi).
\end{equation*}
In general, we obtain
\begin{equation*}
\begin{split}
\V(X_{8;3,1}) &= 4W_8-4V(\chi_{-8}) \\
\V(X_{8;5,1}) &= 4W_8-4V(\chi_{-4}) \\
\V(X_{8;7,1}) &= 4W_8-4V(\chi_8)
\end{split}
\end{equation*}
and
\begin{equation*}
\begin{split}
\V(X_{12;5,1}) &= 4W_{12}-4V(\chi_{-4}) \\
\V(X_{12;7,1}) &= 4W_{12}-4V(\chi_{-3}) \\
\V(X_{12;11,1}) &= 4W_{12}-4V(\chi_{12}).
\end{split}
\end{equation*}
Given the relative sizes of the $V(\chi)$ as listed in
equation~(\ref{Vrelative}), we see that
\begin{equation*}
{\displaystyle
\V(X_{8;5,1}) > \V(X_{8;7,1}) > \V(X_{8;3,1})
 \atop\displaystyle
\V(X_{12;7,1}) > \V(X_{12;5,1}) > \V(X_{12;11,1}).
}
\end{equation*}
This in turn suggests that
\begin{equation*}
{\displaystyle
\Pr(X_{8;3,1}>0) > \Pr(X_{8;7,1}>0) > \Pr(X_{8;5,1}>0)
 \atop\displaystyle
\Pr(X_{12;11,1}>0) > \Pr(X_{12;5,1}>0) > \Pr(X_{12;7,1}>0),
}
\end{equation*}
or equivalently
\begin{equation*}
\delta_{8;3,1} > \delta_{8;7,1} > \delta_{8;5,1} \qquad\text{and}\qquad
\delta_{12;11,1} > \delta_{12;5,1} > \delta_{12;7,1}.
\end{equation*}
This is exactly what is observed in equation~(\ref{2waydeltas}).

We emphasize that although the justification ventured into the analytic realm,
in the end these predictions of the relative sizes of the
$\delta_{q;a,1}$ depended upon only algebraic properties of the various
residue classes $a$ modulo~$q$. To each residue class $a$ was associated a
particular character based on the values of the characters at $a$, and the
conductor of this character is what correlated with the size of
$\delta_{q;a,1}$. This is in the same spirit as Chebyshev's bias: the
sign of $\delta_{q;a,b}-1/2$ was shown by Rubinstein and Sarnak~\cite{RS} to
be determined by whether the residues $a$ and $b$ are squares in the
multiplicative group modulo~$q$.

A similar sort of analysis can also explain the relative sizes of the
densities listed in equation~(\ref{3waydeltas}), for which it is
convenient to define a slightly differently normalized error term for
$\pi(x;q,a)$. When $q=8$ or $q=12$ and $a$ is one of the three
nonsquare residue classes\mod q, we define $\tilde E(x;q,a) = E(x;q,a)
+ E(x;q,1)$; again, investigating the relative sizes of the various
$\pi(x;q,a)$ is the same as investigating the relative sizes of the
$\tilde E(x;q,a)$. For example,
\begin{equation*}
\begin{split}
\tilde E(x;12,5) &= E(x;12,5) + E(x;12,1) \\
&= 2 - \sum\begin{Sb}\chi\mod{12} \\ \chi\ne\chi_0\end{Sb} (\bar\chi(5)+1)
E(x,\chi) + O\big( {1\over\log x} \big)
\\
&= 2 - 2E(x,\chi_{-4}) + O\big( \frac1{\log x} \big),
\end{split}
\end{equation*}
which has the same limiting distribution as the random variable
$2+2X(\chi_{-4})$. In fact, if we define the random variables
\begin{equation*}
\text{
$
\begin{aligned}
\tilde X_{8;3} &= 2 + 2X(\chi_{-8}) \\
\tilde X_{8;5} &= 2 + 2X(\chi_{-4}) \\
\tilde X_{8;7} &= 2 + 2X(\chi_{8})
\end{aligned}
$
\qquad and\qquad
$
\begin{aligned}
\tilde X_{12;5} &= 2 + 2X(\chi_{-4}) \\
\tilde X_{12;7} &= 2 + 2X(\chi_{-3}) \\
\tilde X_{12;11} &= 2 + 2X(\chi_{12}),
\end{aligned}
$}
\end{equation*}
then in each case the distribution of $\tilde E(x;q,a)$ is the same as
that of the random variable~$\tilde X_{q;a}$. Note also that the three
random variables $\tilde X_{8;3}$, $\tilde X_{8;5}$, and $\tilde
X_{8;7}$ are mutually independent due to the hypothesis LI, and the
same is true of $\tilde X_{12;5}$, $\tilde X_{12;7}$, and $\tilde
X_{12;11}$.

If we have three independent random variables each with the same mean,
which one would we expect to take values between those of the other
two most frequently? Our intuition tells us that the random variable
with smallest variance will prefer to stay in the middle, while the
one with largest variance will more frequently be in first or last
place. We easily see that the variances for these random variables are
\begin{equation}
\text{
$
\begin{aligned}
\V(\tilde X_{8;3}) &= 4V(\chi_{-8}) \\
\V(\tilde X_{8;5}) &= 4V(\chi_{-4}) \\
\V(\tilde X_{8;7}) &= 4V(\chi_8)
\end{aligned}
$
\quad and\quad
$
\begin{aligned}
\V(\tilde X_{12;5}) &= 4V(\chi_{-4}) \\
\V(\tilde X_{12;7}) &= 4V(\chi_{-3}) \\
\V(\tilde X_{12;11}) &= 4V(\chi_{12}).
\end{aligned}
$}
\end{equation}
Once again, our knowledge~(\ref{Vrelative}) of the relative sizes of
the quantities $V(\chi)$ tells us that
\begin{equation*}
{\displaystyle
\V(\tilde X_{8;3}) > \V(\tilde X_{8;7}) > \V(\tilde X_{8;5})
 \atop\displaystyle
\V(\tilde X_{12;11}) > \V(\tilde X_{12;5}) > \V(\tilde X_{12;7}).
}
\end{equation*}
Therefore, we expect that of the three prime counting functions
$\pi(x;8,a)$ with $a\in\{3,5,7\}$, the function $\pi(x;8,3)$ spends
more time in first and last place than the other two while the
function $\pi(x;8,5)$ spends the most time in second place; similarly, of the
prime counting functions $\pi(x;12,a)$ with $a\in\{5,7,11\}$, the
function $\pi(x;12,11)$ spends more time in first and last place than
the other two while the function $\pi(x;12,7)$ spends the most time in
second place. All of these predictions match the observed
densities in equation~(\ref{3waydeltas}).

We emphasize how important it was that the trios of random variables
\begin{equation*}
\{\tilde X_{8;3}, \tilde X_{8;5}, \tilde X_{8;7}\}
\quad\hbox{and}\quad \{\tilde X_{12;5}, \tilde X_{12;7}, \tilde
X_{12;11}\}
\end{equation*}
were independent, so that we could draw conclusions about their
relative positions in the three-way race based solely on their
individual variances. We could certainly have normalized the error
terms in an artificial way so that one of the resulting random
variables in a trio equaled zero, for example! But then the other two
random variables would not have been independent, and the correlation
between them would have ruined any chance at such a straightforward
analysis.

We plan to generalize these observations and arguments, as
much as possible, to general moduli $q$ in a future paper. The situation
regarding densities of the form $\delta_{q;a,1}$ for nonsquares $a\mod q$ will
be complicated by the greater complexity of the multiplicative groups to higher
moduli, but we believe that the analysis for the relative sizes of these
two-way densities can be successfully generalized. At the moment, however, the
analysis of the three-way races above relied on the fact that for every
character $\chi\mod q$, at least two of the three values $\chi(a_i)$ were
equal; this is a property that only special triples $\{a_1,a_2,a_3\}\mod q$ can
possess. While these special cases of three-way races to higher moduli can be
treated as above, a new idea will be needed to generalize further.

\section{Densities and equalities}

Since this is a conference proceedings, it seems appropriate to record here
some comments made at the Millennial Conference regarding the
subject of this paper. First, G\'erald Tenenbaum mentioned that the density
\begin{equation}
\delta_{q;a_1,\dots,a_r} = \lim_{x\to\infty}\, {1\over\log x}
\!\!\! \int\limits\begin{Sb}2\le t\le x \\ \pi(t;q,a_1) > \dots >
\pi(t;q,a_r)\end{Sb} \!\!\! {dt\over t},
\label{natdef}
\end{equation}
as defined in equation~(\ref{logdensitydef}) and the preceding lines,
is not the only possible quantity to study when measuring the biases
of the various orderings of the prime counting functions
$\pi(x;q,a_i)$. Indeed, he noted for example that for any real number
$k>-1$, the related density
\begin{equation*}
\delta^{(k)}_{q;a_1,\dots,a_r} = \lim_{x\to\infty}\, {k+1\over(\log x)^{k+1}}
\!\!\! \int\limits\begin{Sb}2\le t\le x \\ \pi(t;q,a_1) > \dots >
\pi(t;q,a_r)\end{Sb} \!\!\!\!\!\! {(\log t)^k\, dt\over t}
\end{equation*}
will also exist in this context.

Surprisingly, it turns out that these densities
$\delta^{(k)}_{q;a_1,\dots,a_r}$ are independent of the parameter $k>-1$. One
can prove this by hand, using the fact that for any function of the form
$f(x) = \alpha x^\beta$ with $\alpha$ and $\beta$ positive, the (natural)
density of those positive real numbers $x$ for which the fractional part of
$f(x)$ lies in an interval $[\gamma,\eta]\subset[0,1]$ is exactly
$\eta-\gamma$. In fact, the lack of dependence on the parameter $k$ is a
consequence of a more general result of Lau~\cite{Lau} regarding
distributions of error terms of number-theoretic functions. So although the
particular definition~(\ref{natdef}) is not canonical, the density values
themselves seems to be natural quantities to consider.

On another topic, Rubinstein and Sarnak~\cite{RS} showed that under
the assumptions GRH and LI, the density of the set of positive real
numbers $x$ such that $\pi(x;q,a)=\pi(x;q,b)$ equals zero (in fact
they prove something rather stronger). Carl Pomerance asked whether
one could prove this particular statement unconditionally. This is an
excellent question, and while it certainly might be possible to
establish unconditionally that $\pi(x;q,a)$ and $\pi(x;q,b)$ are
``almost never'' equal, this author does not know how to do so.

Since we know (conditionally) that the equality $\pi(x;q,a) = \pi(x;q,b)$ has
arbitrarily large solutions, one can ask whether a system of equalities of the
form $\pi(x;q,a_1) = \dots = \pi(x;q,a_r)$ also has arbitrarily large
solutions. A conjecture of the author (see~\cite{FM}), resulting
from an analogy to random walks on lattices, is that the answer is yes when
$r=3$ but no when $r\ge4$. One can refine this conjecture in the following
way: if $\{a_1, \dots, a_r\}$ are mutually incongruent reduced residues\mod q,
then we believe that
\begin{equation}
\liminf_{x\to\infty} \big( \max_{1\le i<j\le r} |\pi(x;q,a_i)-\pi(x;q,a_j)|
\big) = \begin{cases} 0, & \text{if }r\le3, \\
\infty, & \text{if }r\ge4.
\end{cases}
\label{conjecture}
\end{equation}

Of course this raises the issue as to what function $f(x)$ should be chosen
so that for $r\ge4$, the quantity
\begin{equation}
\liminf_{x\to\infty}\, {\max_{1\le i<j\le r}
|\pi(x;q,a_i)-\pi(x;q,a_j)| \over f(x)}
\label{whichproper}
\end{equation}
would be finite and nonzero (and whether the order of magnitude of
this function $f(x)$ depends on~$r$). It follows directly from the
fact that the difference $E(x;q,a)-E(x;q,b)$ possesses a limiting
distribution that for any integers $q,r\ge2$, there exists some
constant $C=C(q,r)>0$ such that the density of those positive real
numbers $x$ with
\begin{equation*}
{|E(x;q,a)-E(x;q,b)|\over\phi(q)} = {|\pi(x;q,a)-\pi(x;q,b)| \over
\sqrt x/\log x} > C
\end{equation*}
is less than $1/r^2$ for any pair $a,b$ of distinct reduced residues
modulo $q$. Thus more than half of the time we must have
\begin{equation*}
{\max_{1\le i<j\le r} |\pi(x;q,a_i)-\pi(x;q,a_j)| \over \sqrt x/\log
x} \le C,
\end{equation*}
since there are only $r(r-1)/2$ terms in the maximum.

This argument shows that the expression in
equation~(\ref{whichproper}) is finite when $f(x)=\sqrt x/\log x$.
However, nothing immediately ensures that the expression is nonzero,
in which case the proper choice of $f(x)$ would be somewhat smaller
than $\sqrt x/\log x$. In any case, we should begin by trying to
establish equation~(\ref{conjecture}) in the first place, perhaps even
in an extreme case such as $r=\phi(q)$.

\end{document}